\documentclass[letterpaper, 10 pt, conference]{ieeeconf}  % Comment this line out if you need a4paper

\usepackage[T1]{fontenc}         
\usepackage{graphicx}
\usepackage{caption}
\usepackage{subcaption} 
\usepackage{float}   
\usepackage{fancybox}		  
\usepackage{listings} % package for a colored Matlab code            
\usepackage{amsfonts}
\usepackage{amsmath} 

\newcommand{\R}{\mathbb{R}}
\newcommand{\N}{\mathbb{N}}
\newcommand{\bi}{\begin{itemize}}
\newcommand{\ei}{\end{itemize}}
\newcommand{\be}{\begin{enumerate}}
\newcommand{\ee}{\end{enumerate}}

\newcommand{\bd}{\begin{description}}
\newcommand{\ed}{\end{description}}
\renewcommand{\i}{\item}

\newcommand{\bqn}{\begin{eqnarray}}
\newcommand{\eqn}{\end{eqnarray}}
\newcommand{\eqnn}{\nonumber\end{eqnarray}}
\newcommand{\eqnl}[1]{\label{#1}\end{eqnarray}}

\newcommand{\nn}{\nonumber\\}

\newcommand{\ba}[1]{\begin{array}{#1}}
\newcommand{\ea}{\end{array}}
\renewcommand{\r}[1]{(\ref{#1})}

\newcommand{\bproof}{\begin{proof}}
\newcommand{\eproof}{\end{proof}}
\newtheorem{Theorem}{\bf Theorem}
\newtheorem{lemma}[Theorem]{\bf Lemma}
\newtheorem{corollary}[Theorem]{\bf Corollary}
\newtheorem{definition}[Theorem]{\bf Definition}
\newtheorem{proposition}[Theorem]{\bf Proposition}
\newtheorem{remark}[Theorem]{\bf Remark}

\newcommand{\bt}{\begin{Theorem}}
\newcommand{\et}{\end{Theorem}}
\newcommand{\bl}{\begin{lemma}}
\newcommand{\el}{\end{lemma}}
\newcommand{\bp}{\begin{proposition}}
\newcommand{\ep}{\end{proposition}}
\newcommand{\bc}{\begin{corollary}}
\newcommand{\ec}{\end{corollary}}
\newcommand{\bdeff}{\begin{definition}}
\newcommand{\edeff}{\end{definition}}
\newcommand{\brem}{\begin{remark}\rm}
\newcommand{\erem}{\end{remark}}
\newcommand{\auth}[1]{{ #1}}
\newcommand{\tit}[1]{{\rm #1}}

\newcommand{\jou}[1]{{\it #1}}
\newcommand{\vol}[1]{{\it #1}}
\newcommand{\pp}[1]{pp.~#1}
\newcommand{\Id}{\mathrm{Id}}

\newcommand{\Pt}[1]{\left( #1 \right)}
\newcommand{\Pg}[1]{\left\{ #1 \right\}}
\newcommand{\Pq}[1]{\left[ #1 \right] }

\newcommand{\schema}[1]{\b{\sc #1}}
\newenvironment{pschema}[1]{\vspace{3mm}
\noindent\begin{Sbox}\begin{minipage}{\columnwidth}\vspace{2mm}\begin{center}{\large \schema{#1}}\vspace{5mm}\\
\begin{minipage}{0.9\textwidth}}{\end{minipage}\end{center}\vspace{2mm}\end{minipage}\end{Sbox}\fbox{\TheSbox}\vspace{3mm}}

\renewcommand{\SS}{\mathbb{S}}
\newcommand{\M}{\mathcal{M}}
\renewcommand{\P}{\mathcal{P}}

\renewcommand{\div}{\mathrm{div}}
\renewcommand{\L}{\mathcal{L}}
\newcommand{\grad}{\mathrm{grad}}

\renewcommand{\b}[1]{{\bf #1}}
\renewcommand{\span}[1]{\mathrm{span}\Pg{#1}}

\newcommand{\Tr}[1]{\mathrm{Tr}\Pt{#1}}

\renewcommand{\Vec}[1]{\mathrm{Vec}\Pt{#1}}

\newcommand{\mmfunz}[5]{
$$ #1 : \left\{ \begin{array}{ccl}
 #2 & \rightarrow & #3 \\
 #4 & \mapsto& #5   \end{array} \right. $$}

\newcommand{\Lip}{\text{Lip}}
\newcommand{\W}{\mathcal{W}}

\IEEEoverridecommandlockouts                              % This command is only needed if 
\title{\LARGE \bf
Control of reaction-diffusion equations on time-evolving manifolds
}

\author{Francesco Rossi$^{1}$,
Nastassia Pouradier Duteil$^{2}$,  
Nir Yakoby$^{3}$, 
Benedetto Piccoli$^{2}$ % <-this % stops a space
\thanks{$^{1}$ Aix Marseille Universit\'e, CNRS, ENSAM, Universit\'e de Toulon, LSIS UMR 7296, 13397, Marseille, France. 
        {\tt\small francesco.rossi@lsis.org}}
\thanks{$^{2}$ Department of Mathematical Sciences and CCIB, Rutgers University - Camden, Camden, NJ.
        {\tt\small piccoli@camden.rutgers.edu, nastassia.pouradierduteil@rutgers.edu}}%
\thanks{$^{3}$ Department of Biology and CCIB, Rutgers University - Camden, Camden, NJ.
        {\tt\small yakoby@camden.rutgers.edu }}
}

\begin{document}

\maketitle
\thispagestyle{empty}
\pagestyle{empty}

%%%%%%%%%%%%%%%%%%%%%%%%%%%%%%%%%%%%%%%%%%%%%%%%%%%%%%%%%%%%%%%%%%%%%%%%%%%%%%%%
\begin{abstract} 
Among the main actors of organism development there are morphogens, 
which are signaling molecules diffusing in the developing organism and acting on cells to produce local responses. Growth is thus determined by the distribution of such signal.
Meanwhile, the diffusion of the signal is itself affected by the changes in shape and size of the organism. 
In other words, there is a complete coupling between the diffusion of the signal and the change of the shapes.\\
In this paper, we introduce a mathematical model to investigate such coupling. The shape is given by a manifold, that varies in time as the result of a deformation given by a transport equation. The signal is represented by a density, diffusing on the manifold via a diffusion equation.
We show the non-commutativity of the transport and diffusion evolution by introducing a new concept of Lie bracket between the diffusion and the transport operator. We also provide numerical simulations showing this phenomenon.
\end{abstract}

%%%%%%%%%%%%%%%%%%%%%%%%%%%%%%%%%%%%%%%%%%%%%%%%%%%%%%%%%%%%%%%%%%%%%%%%%%%%%%%%
\section*{INTRODUCTION}

Modeling growth of living organism attracted the interest of many investigators
both in the field of Developmental Biology and the Applied Mathematics (see \cite{Chauvet, DeV,Jones, Murray, Segelm}).
Developmental biologists have shown that development is primarily induced by morphogens, which act on the organism as  signals by triggering signaling pathways and provoking a response resulting in cell growth or differentiation \cite{Zartman}. 
Several modeling approaches have been explored from the mathematical point of view.
From a microscopic standpoint (see \cite{Jones}), tissues are considered as a collection of cells, and discrete models such as cellular automata are used. We instead adopt a macroscopic standpoint, where the relevant quantity is the density of the signal on a manifold.
As a specific example, in certain fruit flies species such as \textit{Drosophila melanogaster}, a morphogen called Gurken is responsible for initiating the EGFR signalling pathway, resulting in the specification of cells that eventually form structures called dorsal appendages on the \textit{drosphila} eggshell \cite{Goentoro, Zartman}. Interestingly Gurken diffuses in a thin space,
called perivitelline space, which can be modeled by an evolving surface.
This leads naturally to model the growing organism by coupling 
a growing surface with a signal diffusing on  it, see \cite{PRBP-cdc15}.
Because of the biological motivation, this framework was called
 \textit{Developmental Partial Differential Equations}.\\
In this paper we consider a general model, where the boundary of the organism is described by a Riemannian manifold, that evolves with respect to time due to the growth induced by the signal on it. In turn the evolution (for instance, heat diffusion) of the signal on the manifold is affected by the shape of the manifold. Indeed, intrinsic heat diffusion is described by the heat equation with the Laplace-Beltrami operator. Our aim is to investigate the coupling between growth and diffusion.
%For this aim, we built up on the mathematical framework, introduced in \cite{PRBP-cdc15}, and called \textit{Developmental Partial Differential Equations} (DPDE). 
There is a wide literature of studies for PDEs on manifolds,
see for instance \cite{Shatah,Struwe}, or Turing Patterns on evolving manifolds,
see for instance \cite{BM,CHL}. 
However the coupling
of PDE and time-evolving manifold was apparently newly introduced in
\cite{PRBP-cdc15}.\\
As a first step to understand what shapes of the manifold can be attained from an initial configuration, we explore the non-commutativity of the growth
(manifold change in time) and the diffusion operator (on the manifold itself).
A newly defined concept of Lie bracket between the diffusion (2nd order operator) and growth (1st order operator) is able to capture such non-commutativity
and thus provide new shapes towards which the manifold may evolve.
As in classical geometric control theory \cite{agrabook,bressan-piccoli,sontag}, the concept of Lie bracket may indeed enclose all the needed information to capture the controlled dynamics.
Moreover, such bracket can be understood as a new available direction for the growth of the organism.
% by slowly perturbing the diffusion and growth.

The paper is structured as follows. We begin by introducing the general model, or Developmental Partial Differential Equation describing the coupling of growth and diffusion on a Riemannian manifold. We then prove existence and uniqueness of the solution to the DPDE by introducing a numerical scheme that discretizes time and solves diffusion and growth independently on each time interval. We prove that the limit of the scheme is the solution to the DPDE. 
We then use the scheme to define a new kind of Lie bracket between the diffusion and the growth operators. By computing the bracket explicitly, we show that it is not zero.
Numerical simulations confirm the analytical computation of the bracket.

\section{DESCRIPTION OF MORPHOGENESIS}

In this section, we describe a simplified model for morphogenesis, i.e. for the development of the shape of a living body.
The shape of an organism is described by its boundary, represented by a time-varying manifold $\M_t$ embedded in an Euclidean space $\R^d$ with the dimension $d$ being fixed (naturally $d=3$ in real examples). On such manifold, a growth signal is represented by a probability measure $\mu_t\in \P(\M_t)$. Here $\P(\M_t)$ is the space of probability measures on $\M_t$, endowed with the Wasserstein distance $W_p$ (see more details in Section \ref{s-wass}).  
Using the embedding of $\M_t$ into the ambient space $\R^d$, we can consider $\mu_t$ as a probability measure on $\R^d$.\\
The organism development is determined by a growth vector field $v[\mu_t]$ given by the current shape of the organism and by the signal, with $v[\cdot]:\P_c(\R^d)\rightarrow \Lip(\R^d)\cup\mathcal{L}^\infty(\R^d)$. 
The signal $\mu_t$ on $\M_t$ diffuses following the heat equation intrinsically defined on $\M_t$. Indeed, since $\M_t$ inherits the Riemannian structure of the ambient space $\R^d$,  we can define an intrinsic Laplacian, called the Laplace-Beltrami operator. We denote with $\Delta_t$ the Laplace-Beltrami operator on $\M_t$.\\
These two phenomena (growth and diffusion) can be summarized by describing the evolution of the signal $\mu_t$ via the following transport-diffusion Partial Differential Equation :
\bqn
\partial_t\mu_t+\nabla\cdot(v\Pq{\mu_t} \mu_t)=\Delta_t \mu_t,
\eqnl{e-main}
where the manifold $\M_t$ is the support of $\mu_t$ at each time. Since $\mu_t$ are measures in $\R^d$, such equation needs to be interpreted in the weak sense, i.e. for all $f\in C^\infty(\R^d)$ it holds
 \begin{equation}\label{eq:weak}
 \partial_t \int_{\R^d} f d\mu_t - \int_{\R^d} (\nabla f \cdot \nabla v[\mu_t])d\mu_t = \int_{\R^d} \Delta_t f \; d\mu_t.
 \end{equation}
We will provide in Section \ref{s-ex} existence and uniqueness results for such equation. 
\begin{remark}
Notice that the developed theory can be adapted to include reaction
terms of the type $h[\mu_t]$ on the right-hand side of (\ref{e-main})
by using the generalized Wasserstein distance, see \cite{piccoli-rossi}.
\end{remark}

 \subsection{Transport equation and Wasserstein distance}\label{s-wass}
 In this section, we recall the definition of the Wasserstein distance and its connection with nonlinear transport equation, i.e. equation \r{e-main} with no diffusion. 
Let us first recall the definition of the Wasserstein distance (see \cite{villani}).
Recall that, for every probability measure $\mu$ and measurable map $\phi$,
the push-forward $\phi\#\mu$ is defined by $\phi\#\mu(A)=\mu(\phi^{-1}(A)).$
\begin{definition} Fix $p\geq 1$. Given two probability measures $\mu$ and $\nu$ in $\R^d$, the $p$-Wasserstein distance between $\mu$ and $\nu$ is given by:
\bqn
%\W(\mu,\nu) := \sup\{ \int_{\R^d}fd(\mu-\nu)\ \; | f \in C^\infty(\R^d), \Lip(f)\leq 1\}.
\W_p(\mu,\nu) := \min\limits_{\pi\in\Pi(\mu,\nu)}\Pt{\int_{\R^d\times\R^d} |x-y|^p d\pi(x,y)}^{1/p}
\eqnn
where $\Pi(\mu,\nu)$ is the set of transference plans from $\mu$ to $\nu$, i.e. of the probability measures on $\R^d\times\R^d$ with marginals $\mu,\nu$, respectively. In other words $P_x\#\pi=\mu$ and $P_y\#\pi=\nu$
(where $P_x$, respectively $P_y$ denote the projection on the first,
respectively second, component of $(x,y)$.)
\end{definition} 
The transference plans in $\Pi(\mu,\nu)$ can be seen as methods to transport 
$\mu$ to $\nu$ and the term $\int_{\R^d\times\R^d} |x-y|^p d\pi(x,y)$ can be interpreted as a cost (as $p$-power of the distance) to move 
the mass of $\mu$ onto the mass of $\nu$ via the plan $\pi$. 
Hence, the Wasserstein distance is the minimal cost to move one mass over the other.  For a complete introduction to the topic of Wasserstein distances 
we refer the reader to \cite{villani}.\\
Let us now consider the Cauchy problem
 \begin{equation}\label{e-trasporto} 
 \begin{cases}
\partial_t\mu_t+\nabla\cdot(v\Pq{\mu_t} \mu_t)=0,\\
\mu(t=0)=\mu_0.
\end{cases}
 \end{equation}
 We assume that $v$ is a uniformly Lipschitz operator with respect to the Wasserstein distance on $\P(\R^d)$ and  the Euclidean distance in $\R^d$, i.e. that there exists a constant $L$ such that 
\bqn
\|v[\mu](t,x)-v[\nu](t,y)\|\leq LW_p(\mu,\nu)+L\|x-y\|
\eqnl{e-vlip}
for all $t\in \R$, $\mu,\nu\in \P(\R^d)$ and $x,y\in\R^d\times\R^d$.
We have the following key result (\cite{pedestrian}):
\bt
Let $v$ satisfy \r{e-vlip}, then there exists a unique solution to \r{e-trasporto}.\et

\subsection{Existence of a solution to \r{e-main}} \label{s-ex}
In this section, we prove existence of a solution for \r{e-main}, by providing a numerical scheme approximating such solution.
Fix a final time $T\in\R$ and an initial measure $\mu_0$. For a given discretization parameter $n\in\N$, we define a sequence of curves $(\mu_s^n)$ via the following scheme:
\begin{pschema}{Scheme $\SS$}
Define $\tau_n=t_n:=2^{-n}T$.
Let $\mu^n(0):=\mu_0$.
On the nodes $lt_n$ (with $l\in \{0,...,2^n-1\}$) we define $\mu^n((l+1)t_n)$ from $\mu^n(lt_n)$ as follows:
\be
\i Let $\phi_{lt_n}^{t_n}$ be the flow of $v[\mu^n(lt_n)]$ and consider
$\phi_{lt_n}^{t_n}\#\mu^n(lt_n)$, i.e. the push-forward of $\mu^n(lt_n)$ via the flow $\phi_{lt_n}^{t_n}$, that is a measure on $\M_{(l+1)t_n}$.
\i Define $\mu^n((l+1)t_n)=e^{\Delta_{lt_n}\tau_n}(\phi_{lt_n}^{t_n}\#\mu^n(lt_n)),$ i.e the solution of the heat equation on $\M_{(l+1)t_n}$ with initial data $\phi_{lt_n}^{t_n}\#\mu^n(lt_n)$ at time $\tau_n$.
\ee
In between nodes, for $t\in (0,t_n)$ we define:\\
$\mu^n(lt_n+t) = e^{\Delta_{lt_n+t}t}(\phi^t_{lt_n}\#\mu^n(lt_n))$.
\end{pschema}
In the definition of $\SS$, we distinguish $t_n$ and $\tau_n$ for better description and approximation of the two phenomena of deformation and heat diffusion.
We now prove existence of a solution to \r{e-main} with the following lemma. 
\begin{lemma}
%Let $s\in [0,T]$. Then 
There exists a subsequence of $(\mu^n)$ converging to a measure $\mu^*$, 
providing a solution of \r{e-main}.
%:=\lim_{n\to \infty} \mu^n$ . 
%Moreover, $\mu^*$ is a continuous curve in $\Pr$, satisfying $\mu^*_0=\mu_0$.
\end{lemma}
\bproof 
We will prove that each $\mu^n$ is H\"older as function of time
with values in the space of probability measures (endowed
with the Wasserstein distance), in order to use the Arzel\`a --Ascoli theorem. 
The proof will be performed for $p=2$ thus we drop the subscript $p$
and simply write $\W$.
Let $l\in \{ 0,...,2^n-1\}$, $t\in (0,t_n)$ and define $\sigma:=\phi^t_{lt_n}\#\mu^n(lt_n)$. Then by the triangular inequality,
\bqn
\begin{split}
& \W(\mu^n(l t_n+t),\mu^n(l t_n)) = \W( e^{\Delta_{lt_n+t}t}\sigma, \mu^n(l t_n))\\
\leq & \; \W( e^{\Delta_{lt_n+t}t}\sigma, \sigma) + \W(\sigma,\mu^n(l t_n)).
\end{split}
\eqnl{e-conv}
For the first term, we use the evolution variational inequality given in \cite{Erbar}:
\bqn
\begin{split}
 & \frac{d}{dt}\frac12 \W^2( e^{\Delta_{lt_n+t}t}\sigma, \sigma) \\
\leq & \; H(\sigma)-H(e^{\Delta_{lt_n+t}t}\sigma)  -\frac{K}{2}\W^2( e^{\Delta_{lt_n+t}t}\sigma, \sigma)
%\leq & \; H(e^{\Delta_{lt_n+t}t}\sigma) -H(\sigma)
\end{split}
\eqnn
where $H(\rho)$ denotes the relative entropy of $\rho$ and $K$ the lower bound of the Ricci curvature of $\M$. By recalling that the heat equation is the gradient flow for the relative entropy $H$, it holds $H(e^{\Delta_{lt_n+t}t}\sigma)\leq H(\sigma)$. Also observe that the relative entropy is bounded on a compact manifold, see e.g. \cite[Lemma 4.1]{Erbar}. Finally, observe that the Ricci curvature $K$ is bounded from below on a compact manifold too. Hence we obtain :
\bqn
 \frac{d}{dt} \W^2( e^{\Delta_{lt_n+t}t}\sigma, \sigma)
\leq C -K \W^2( e^{\Delta_{lt_n+t}t}\sigma,\sigma)
\eqn
where $C$ is independent of $n$ and, by Gronwall's inequality: 
\bqn
\W^2( e^{\Delta_{lt_n+t}t}\sigma, \sigma)\leq -\frac{C}{K}(e^{-Kt}-1) = C t +o(t) \leq 2Ct
\eqnn
for $t$ small enough. Finally, for some $C>0$ it holds:
\bqn
\W( e^{\Delta_{lt_n+t}t}\sigma, \sigma)\leq C\sqrt{t}.
\eqnn
The second term of \r{e-conv}  was estimated in \cite{pedestrian}: $\W(\phi^t_{lt_n}\#\mu^n(lt_n),\mu^n(l t_n))\leq L t$ where $L$ is the Lipschitz constant of $V[\mu_s]$. Notice that $L$ does not depend on $\mu^n$ nor on $t$. Summing the two terms we obtain: 
\bqn
\W(\mu^n(l t_n+t),\mu^n(l t_n))\leq Lt+ C\sqrt{t}.
\eqnn
Iteratively and by the triangular inequality, we get:
$\W(\mu^n(t),\mu^n(s))\leq L|t-s|+ C\sqrt{|t-s|}$ for any $s,t\in[0,T]$.
Hence the sequence $(\mu^n)$ satisfies a uniform H\"older condition of order $\frac{1}{2}$. The sequence is also equibounded, since $\mu_0$ is fixed. Then, according to the Arzel\`a --Ascoli theorem, there exists a subsequence of $(\mu^n)$ that converges uniformly to a curve $\mu^*$.
Using the same methods as \cite{pedestrian}[Sec 3.3], one can prove that $\mu^*$ is a solution to \r{e-main}.
\eproof

%\brem The construction above gives a solution to \r{e-main}. Nevertheless, uniqueness is not guaranteed, since the scheme could have more than one limit. Moreover, it could be possible to have a solution not given by the limit of such scheme. 
%\erem

\section{DEFINITION OF LIE BRACKET}

% In the following section, we define a bracket between the diffusion operator and the growth vector field.
Here onward, for simplicity we assume that the growth vector field $v$ does not depend on $\mu_t$. This is a suitable approximation since the Lie bracket
we are going to define is a local object (as the original Lie bracket).

\subsection{Reduction to a time-varying Riemannian structure on $\M_0$}
\label{s-trasf}
In this section, we transform the problem \r{e-main} on a time-varying manifold $\M_t$ into a problem defined on the fixed manifold $\M_0$ with a time-varying Riemannian structure. We use this transformation to prove uniqueness of the solution to \r{e-main} with $v$ not depending on $\mu_t$.\\
Given $v$ not depending on $\mu_t$, the definition of $\M_t$ is given by $\M_t=\phi^t(\M_0)$, where $\phi^t$ is the flow of $v$.
Endow $\M_t$ with the Riemannian structure given by its embedding in $\R^d$. For each time $t$, the flow $\phi^t$ is a diffeomorphism between $\M_0$ and $\M_t$, hence we can endow $M_0$ with a Riemannian structure induced by the one on $\M_t$. Applying this technique at each time, we have defined a time-varying Riemannian structure $<.,.>_t$ on the fixed manifold $\M_0$. We denote with $\Delta'_t$ the corresponding Laplace-Beltrami operator on $\M_0$.
We are now ready to prove uniqueness of solution to \r{e-main}.
\bt
Let $v$ be a Lipschitz vector field on $\R^d$, independent of $\mu$. Then, there exists a unique solution to \r{e-main}.
\et
\bproof Let $\mu_t$ any solution to \r{e-main}, define the measure $\nu_t:=\phi^{-t}\#\mu_t$, i.e. the pull-back of $\mu_t$ via the flow $\phi^t$ generated by $v$. The transformation defined above permits to  prove that $\nu_t$ is a measure on $\M_0$ and it satisfies the following heat equation on $\M_0$:
$$\partial\nu_t=\Delta'_t\nu_t.$$
% See more details in \cite{luca}. 
Observe now that such equation admits a unique solution, see e.g. \cite{guenther}. Moreover, uniqueness of $\nu_t$ implies uniqueness of $\mu_t$.
\eproof

\subsection{First-order Taylor expansion of the Laplace-Beltrami operator}

In this section, we describe the evolution of the Laplace-Beltrami operator $\Delta_t$ by computing its first-order Taylor expansion at time $t=0$.
By applying the transformation described in Section \ref{s-trasf}, we consider the Laplace-Beltrami operator as being defined on the fixed manifold $\M_0$, with time-varying Riemannian structure $<.,.>_t$ on it. For simplicity of notation, we denote $\M:=\M_0$.\\
We now compute the Riemannian structure $<.,.>_t$ on $\M$. The definition of the push-forward implies
\bqn
<w_1,w_2>_t=<\phi_t \# w_1,\phi_t \# w_2>_E,
\eqnl{e-metric}
where $<.,.>_E$ is the standard Riemannian metric in $\R^d$
(i.e. the standard scalar product.)
% (or eventually any other Riemannian metric in $\R^d$). 
It holds
\bqn
\phi_t \# w=w+t Jv\cdot w+o(t),
\eqnn
where $J$ is the Jacobian with respect to the Euclidean structure of $\R^d$ and $\cdot$ represents the linear action of the linear operator $Jv$ on $w$. This implies
\bqn
\begin{split}
<w_1,w_2>_t=& <w_1,w_2>_E+t (<Jv\cdot w_1,w_2>_E \\
& +<Jv\cdot w_2,w_1>_E)+o(t).
\end{split}
\eqnn
Since vectors $w_1,w_2$ belong to $T_x\M$, we will denote with $J_{\M}v$ the restriction of $Jv$ to $T_x\M$ by projection, i.e.
\mmfunz{J_{\M}v}{T_x\M}{T_x\M}{w}{(Jv\cdot w)_{\M},}
where $z_{\M}$ is the component of the vector $z\in T_x\R^d$ on the subspace $T_x\M$. Observe that we are using here the Riemannian structure of $\R^d$ to define projections.\\
We now compute the Laplace-Beltrami operator $\Delta'_t$ intrinsically defined on the manifold $\M$ with the Riemannian structure $<.,.>_t$. We are interested in describing such operator as a function of $t$. Recalling that the Laplace-Beltrami operator is the divergence of the gradient, we aim at describing $\div^t$ and $\grad^t$ as a function of time. In particular, we aim at computing first-order development of such operators with respect to time.
We first study the gradient $\grad^t f$ for a function $f\in C^\infty(\M)$, via its intrinsic definition. For all $w\in T_x \M$ it holds
\bqn
<\grad^t(f),w>_t=\L_w f.
\eqnn
In particular this identity holds at time $t=0$ holds for $\grad^0$. 
Writing $\grad^t f=\grad^0 f+t B_1$, for a vector field $B_1$ to be found, 
we get
\bqn
\begin{split}
& <\grad^0 (f)+tB_1+t Jv\cdot \grad^0 (f)+o(t),\\
& \quad w+t\, Jv\cdot w+o(t)>_E \\
= & \L_w f t (<B_1,w>_E+<Jv\cdot \grad^0 (f), \\
& w>_E+<\grad^0 f, Jv\cdot w>_E)+o(t)=0.
\end{split}
\eqnn
We then have $B_1=-(Jv\cdot \grad^0(f))_{\M}-B$ where $(Jv\cdot \grad^0(f))_{\M}$ is the component of $Jv\cdot \grad^0(f)$ on the tangent space of $\M$, and $B(f,v)$ is intrinsically defined by the following rule: for all $w\in T_x\M$ it holds
\bqn
<B(f,v),w>_E=<\grad^0 (f),Jv\cdot w>_E.
\eqnl{e-B}
Summing up, we have
\bqn
\grad^t(f)=\grad^0(f)-t (Jv\cdot \grad^0(f))_\M-t B(f,v)+o(t)
%\eqnl{e-grad}
\eqnn
with $B(f,v)$ defined by \r{e-B}.

\renewcommand{\vol}{\mathrm{vol}}

We now study the divergence $\div^t(X)$ for a vector field $X\in\Vec{\M}$. 
Denoting by $\vol_t$ the volume form of the Riemannian manifold, it holds
\bqn
\div^t(X)\vol_t=\L_X \vol_t.
\eqnl{e-div1}
Observe that $\vol_t =\sqrt{|g^t|}\,dX_1\wedge dX_2 \wedge \ldots\wedge dX_m$ for any base $X_1,\ldots,X_m$ of the Riemannian manifold $(M,<.,.>_t)$. We choose an orthonormal basis for $(M,<.,.>_0)$ and study the evolution of $\vol_t$. Since $g^0=\Id$ and $|\Id+tA|=1+t\Tr{A}+o(t)$, we have
\bqn
|g^t|= 1+2t\sum_{i=1}^m <Jv\cdot X_i,X_i>_E=1+2t\Tr{J v}_{\M},
\eqnn
where the operator $Jv$ is restricted to the tangent space of $\M$. 
Then $\sqrt{|g^t|}=1+t\Tr{J v}_{\M}$. Writing $\div^t(X)=\div^0 X+t f+o(t)$ for a function $f$ to be found, we get
\bqn
\begin{split}
& ((\div^0(X)+t f)(1+t \Tr{J v}_{\M})+o(t))\vol_0 \\
=&\Pt{\L_X (1+t\Tr{J v}_{\M})} \vol_0 +(1+t\Tr{J v}_{\M})\L_X \vol_0.
\end{split}
\eqnn
% =&\Pt{\L_X (1+t\Tr{J v}_{\M})} \vol_0 +(1+t\Tr{J v}_{\M})\div^0(X) \vol_0
From \r{e-div1} applied for $t=0$ we have $\L_X \vol_0=\div^0(X) \vol_0$, thus
\bqn
\div^t(X)=\div^0(X)+t \L_X\Tr{J v}_{\M}.
\eqnl{e-div}
Observe that this formula is intrinsic, since the trace of the linear operator $J v$ does not depend on the chosen orthonormal frame.

We now compute the Laplace-Beltrami operator $\Delta'_t$. Since $\Delta'_t f=\div^t(\grad^t(f))$ by definition, and observing that it holds $L_{\grad^0(f)}\Tr{J v}_\M=<\grad^0(f),\grad^0(\Tr{J v}_\M)>_E$, we have
\bqn
\begin{split}
& \Delta'_t(f)=\Delta_0(f)+t(<\grad^0(f),\grad^0(\Tr{J v}_\M)>_E \\
& -\div^0(B(f,v)+(Jv\cdot \grad^0(f))_\M))+o(t).
\end{split}
\eqnl{e-devel}

\subsection{Non commutativity of the heat and growth evolutions}

In this section, we show non-commutativity of the growth and diffusion terms in the dynamics \r{e-main}, which describes morphogenesis.
%Let us assume that one is allowed to vary slightly the values of $t_n,\tau_n$ (that are no more identical) in the scheme $\SS$. Then, we will prove then the final value of the scheme $\SS$ will be different than the solution of \r{e-main}. In terms of reachability: authorising a slight variations of $t_n,\tau_n$ implies a larger set of attainable final configurations. The goal of this section is to prove that the set of attainable configurations is larger, and to describe such set.
For sake of clarity, let us recall how non-commutativity works in the finite
dimensional case, and, more specifically for switching systems.
%Recall that, for finite dimensional systems, non-comutativity is studied for the so-called switching systems. 
Consider two vector fields $X_0,X_1$ and the system
\bqn
\begin{cases}
\dot x=X_{u(t)}\\
x(0)=x_0.
\end{cases}
\eqnn
where $u\in \{0,1\}$. 
For every measurable switching function $u:[0,T]\to \Pg{0,1}$ 
the solution at time $T$, denoted by $x(T,u)$, is unique. 
The set of points reachable with these trajectories is not limited
to the directions given by $X_0$ and $X_1$.
%Nevertheless, if one chooses another switching function $v:[0,T]\to \Pg{0,1}$, one has in general $x(T,u)\neq x(T,v)$. Then, it is of interest to study the set of all possible final states, at least for small $T$. 
A classical result in control theory, the Orbit theorem, (roughly) states that the set of attainable configurations is related to the Lie bracket $[X_0,X_1]$
(and to other higher order brackets), and in particular that one can choose good switching functions to drive the system along a direction arbitrarily close to the vector field $[X_0,X_1]$. See e.g. \cite{bressan-piccoli,sontag,agrabook}.
Lie brackets are a powerful tool for finite-dimensional control systems. In infinite dimension, they are less easy to define. However, in some cases, they can still be used to study the controllability of PDEs, as shown in \cite{Coron}.\\
For this reason, we study in this article the bracket between the ``heat vector field'' and the ``transport vector field''. Indeed, one can consider the solution of an heat equation as a continuous (and even differentiable) curve in $P_2(X)$ endowed with the Wasserstein distance. The, the time derivative of this curve in a point $\mu_t$ (that is clearly the Laplacian $\Delta_t\mu_t$) can be considered as a vector field, that we call the {\bf heat vector field}. Similarly, we define the {\bf transport vector field} as the derivative of the solution of the transport equation in a point.\\
By borrowing the notation from Lie brackets of vector fields, we define
\bqn
[\Delta,v]\mu:=\lim_{t,\tau\to 0}\frac{\phi_{-t}\#\Pt{e^{\tau\Delta_t}(\phi_t\#\mu)}-e^{\tau\Delta_0}\mu}{t\tau},
\eqnl{e-bracket1}
where $\phi_t$ is flow generated by the vector field $v$, and $e^{\tau \Delta_t}$ is the semigroup generated by $\Delta_t$ at time $\tau$.
%Observe that, for any measure $\mu\in \M(\M_0)$, it holds $\Phi_t\#\mu\in \M_0(\R^{2d}(\M_t)$. If one chooses the coordinates on the manifold $\M_t$ induced by the diffeomorphism $\Phi_t$, then it holds $\Phi_t\#\mu=\mu$ in the sense that their expression with respect to coordinates is the same. Similarly, it holds $\Phi_{-t}\#\mu=\mu$ for any measure $\mu\in\Mu_0(\R^{2d}(\M_t)$. 
Then, for any test function $f\in C^\infty_c(M)$, one can write \r{e-bracket1} as follows:
\bqn
&([\Delta,v]\mu)(f) =
\eqnn
\bqn
&=&\lim_{t,\tau\to 0}\int_\M \frac{f}{t\tau}\,d \Pt{\phi_{-t}\#\Pt{e^{\tau\Delta_t}(\phi_t\#\mu)}-e^{\tau\Delta^0}\mu}\nn
&=&\lim_{t,\tau\to 0}\int_\M\frac{f}{t\tau}\,d \Pt{e^{\tau\Delta'_t}\mu-e^{\tau\Delta_0}\mu}\nn
&=&\lim_{t\to 0}\int_\M \frac{f}{t}\,d \Pt{\Delta'_t\mu-\Delta_0\mu}=\lim_{t\to 0}\int_\M (\Delta'_t-\Delta_0)\frac{f}{t}\,d \mu \nn
&=&\int_\M <\grad^0(f),\grad^0(\Tr{J v}_\M)>_E \nn
& &\quad -\div^0(B(f,v) +(Jv\cdot \grad^0(f))_\M)\,d\mu,
\eqnl{e-brackint}
where we used $\int f\,d \Delta \mu=\int \Delta f\,d\mu$,
which is the definition of the Laplace-Beltrami operator as an operator on the space of measures. Then, \r{e-brackint} is the intrinsic formula for the bracket \r{e-bracket1}, indeed the operator $B$ is intrinsically defined by \r{e-B}.

\section{SIMULATIONS}
\begin{figure}[h!]
        \begin{center}
                \includegraphics[trim=2cm 1cm 0cm 0.0cm, clip=true, scale=0.5]{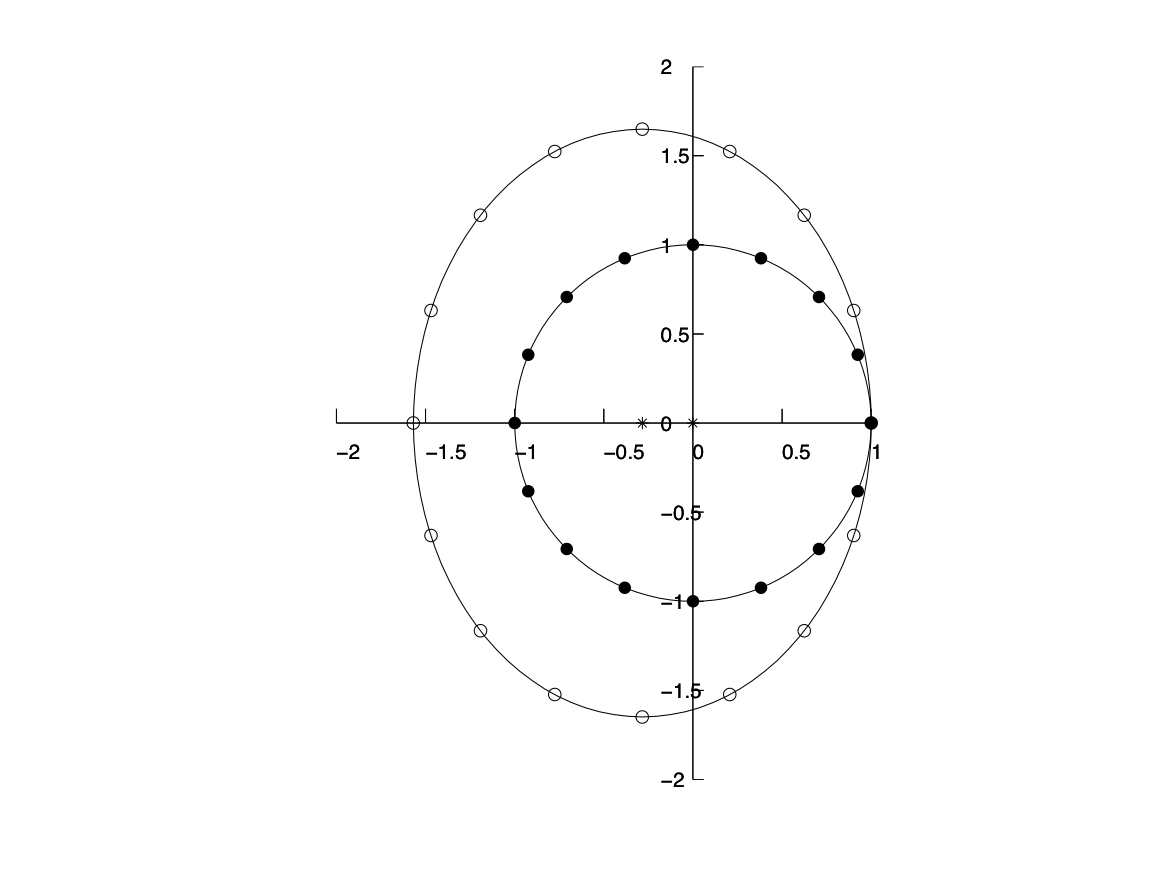} 
        \caption{Transport of the unit sphere (black dots) by the vector field $v(x,y):=(x-1,2y)$. At $t=0.25$, the resulting ellipse (white dots) is centered at $(1-e^{0.25},0)$. }\label{fig:S1transport}
        \end{center}
\end{figure}

\begin{figure}[h!]
        \begin{center}
       % \begin{subfigure}[b]{0.2\textwidth}
                \includegraphics[trim=0cm 0cm 0cm 0cm, clip=true, scale=0.4]{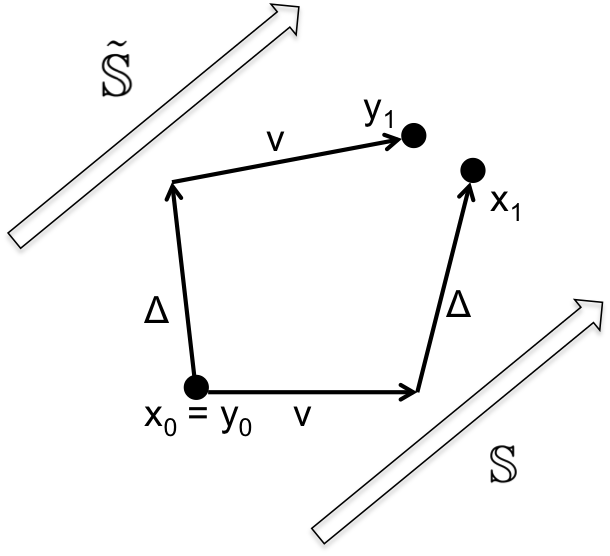}
        \caption{One iteration of the schemes $\SS$ and $\tilde{\SS}$ starting from the same point $x_0=y_0$.}\label{fig:Schemes}
        \end{center}
\end{figure}

\begin{figure}[h!]
        \begin{center}
                \includegraphics[trim=0cm 0.3cm 0cm 0.5cm, clip=true, scale=0.4]{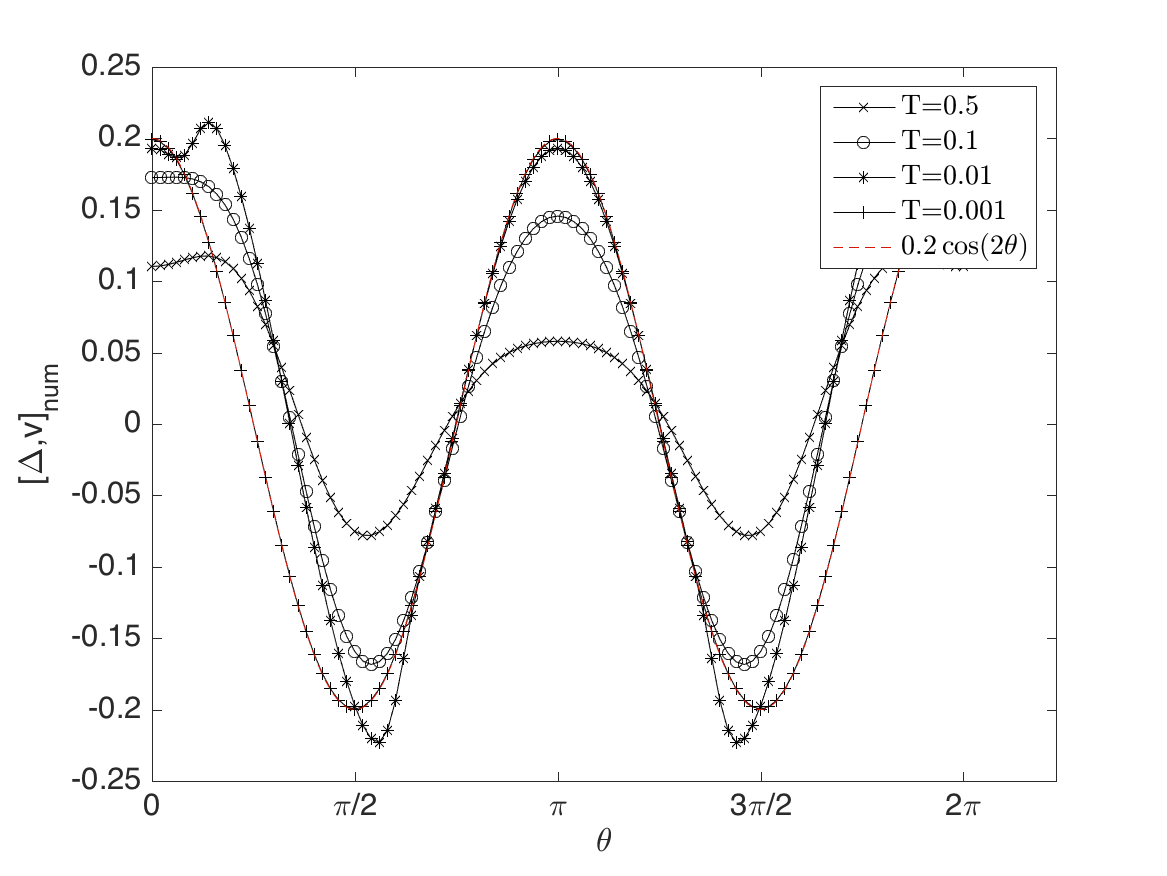}
        \caption{Convergence of the bracket to the theoretical one for the initial signal $\mu_0(\theta)=0.1d\theta$ .}\label{fig:BracketnumConstantSignal}
        \end{center}
\end{figure}

\begin{figure}[h!]
        \begin{center}
                \includegraphics[trim=0cm 0.3cm 0cm 0.5cm, clip=true, scale=0.4]{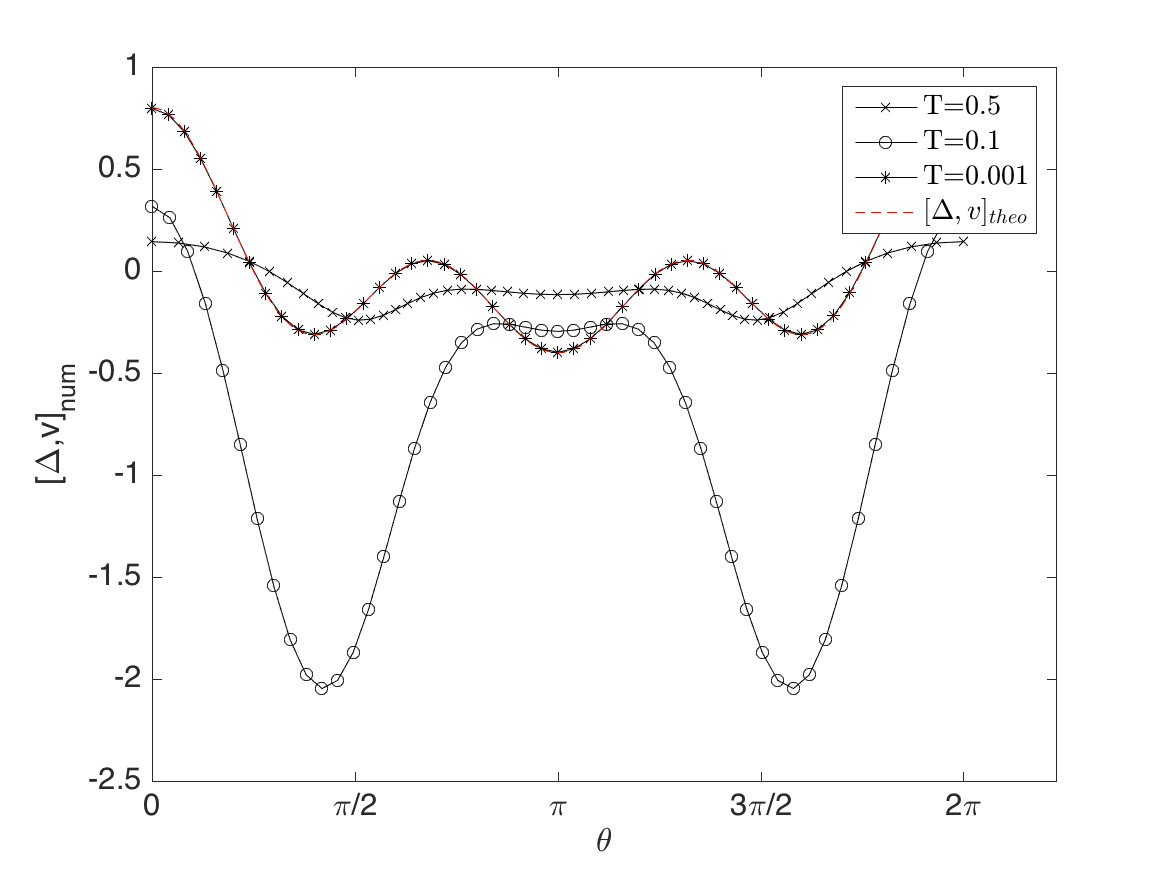}
        \caption{Convergence of the bracket for the initial signal $\mu_0(\theta)=0.1(\cos(\theta)+1)d\theta$.}\label{fig:BracketnumCosineSignal}
        \end{center}
\end{figure}

We now compute the bracket $[\Delta,v]$ for two examples on a time-varying manifold $\M_t$ with $v$ independent on $\mu$. 
We will also compare the analytic computation using formula \r{e-brackint} and the numerical simulations with the scheme $\SS$ defined in Section \ref{s-ex}.\\
We consider the unit circle $\M_0=S^1$ in $\R^2$ parametrized by an angle $\theta$ as the initial manifold, and the vector field $v=(x-1,2y)$. It is easy to verify that at time $t$ the unit circle is transported to an ellipse of equation: $\left( \frac{x-x_c}{e^{ t}}\right)^2+\left( \frac{y}{e^{2 t}}\right)^2=1$ where $x_c=1-e^{ t}$ (see Fig. \ref{fig:S1transport}). 
We consider the Euclidean metric on $\R^2$, i.e. Riemannian structure given by the orthonormal frame $\partial_x,\partial_y$ at each point. The corresponding Riemannian structure on $M_t$ is given by $\partial_\theta=-y\partial_x+x\partial_y$. This implies
\bqn
(J_Mv)\cdot \delta_\theta=\Pt{\ba{c}-y\\2x\ea}_M=(1+x^2)\partial_\theta,
\eqnn
thus $\Tr{J_Mv}=1+\cos^2$. Since the initial data is the Riemannian volume form, the divergence theorem implies
\bqn
([\Delta,v]\mu)(f)&=\int_M <\frac{\partial f}{\partial \theta} \partial_\theta,\frac{\partial (1+\cos^2(\theta))}{\partial \theta}\partial_\theta>_E\,d\theta \nn
=&\int_M -f\frac{\partial^2 \cos^2(\theta)}{\partial \theta^2}\,d\theta=\int_M 2 f\cos(2\theta)\,d\theta.
\eqnn
As a first exemple, we consider an initial constant signal $\mu_0(\theta)=\frac{1}{2\pi}$. We then have $[\Delta,v]\mu_0=2\cos(2\theta)\mu_0$.
As a second exemple, for a more complicated initial data $\mu_0=(1+\cos(\theta))d\theta$, the second term in \r{e-brackint} is no longer 0. First, notice that:
\bqn
\begin{split}
& <\grad^0(f),\grad^0(\Tr{J v}_M)> \\
= & <\frac{\partial f}{\partial \theta} \partial_\theta,\frac{\partial (1+\cos^2(\theta))}{\partial \theta}\partial_\theta> 
 =  \frac{\partial f}{\partial \theta}(-2\cos\theta \sin\theta).
\end{split}
\eqn
Secondly, we calculate the term $B(f,v)$ knowing that $\langle B(f,v),w)\rangle_E = \langle \grad^0(f),Jv\cdot w\rangle_E.$
Taking a vector $w = (w_x,w_y)^T$, we write:
\bqn
\begin{split}
& B_x w_x + B_y w_y = \langle \grad^0(f),Jv\cdot w\rangle_E \\
& = \langle 
\begin{pmatrix}
  \partial_x f \\
\partial_y f 
 \end{pmatrix},
 \begin{pmatrix}
  w_x  \\
2w_y  
 \end{pmatrix} \rangle_E 
  = \partial_x f w_x +2 \partial_y f w_y.
  \end{split}
\eqn
Hence we have : 
$$
B_\theta = \langle \begin{pmatrix}
  \partial_x f \\
2 \partial_y f 
 \end{pmatrix}, 
 \begin{pmatrix}
  -y \\
x f 
 \end{pmatrix}\rangle  
 = (1+\cos^2\theta)\partial_\theta f.
$$
Similarly, 
$$
(Jv\cdot \grad_0(f))_M = \langle
 \begin{pmatrix}
  \partial_x f \\
2 \partial_y f 
 \end{pmatrix}, \partial_\theta\rangle\partial_\theta  = (1+\cos^2\theta)\partial_\theta f \; \partial_\theta.
$$
So
\bqn
\begin{split}
&\div^0 (B(f,v)+Jv\cdot\grad^0(f))_M \\
= & 2\partial_\theta ((1+\cos^2\theta)\partial_\theta f) \\
 = & 2(1+\cos^2\theta) \partial_\theta^2 f -4\cos\theta\sin\theta \partial_\theta f. 
\end{split}
\eqnn
%We calculate separately the two terms in Equation \eqref{e-brackint} :
%\bqn
%%\begin{split}
%& & \int_M <\grad^0(f),\grad^0(\Tr{J v}_M)> (1+\cos\theta) d\theta= \nn
%%& = \int_M \frac{\partial f}{\partial \theta} (-2\cos\theta \sin\theta) (1+\cos\theta) d\theta \\ 
% %& = \int_M 2 f \frac{\partial}{\partial \theta} ((\cos\theta \sin\theta) (1+\cos\theta)) d\theta \\
%&& = \int_M  f (6\cos^3\theta + 4 \cos^2\theta -4 \cos\theta -2) d\theta.\nn
%& & \int_M - \div^0((B(f,v) + Jv\cdot \grad^0(f))_M) (1+\cos\theta)d\theta \nn
%%& = - 2 \int_M ( (1+\cos^2\theta)\frac{\partial^2f}{\partial\theta^2} - 2 \cos\theta\sin\theta\frac{\partial f}{\partial \theta})\cos\theta d\theta \\ 
% %& = 2 \int_M [ \frac{\partial f}{\partial \theta} (-\sin\theta-3\cos^2\theta\sin\theta) \\
% %& \quad -2f (-2\cos\theta\sin^2\theta+\cos^3\theta)]d\theta \\
%&& = 2 \int_M f (-\cos\theta+3\cos^3\theta) d\theta.
%%\end{split}
%\eqnn
%Summing the two terms we get:
Then, a direct computation shows
\bqn
([\Delta,v]\mu)(f)=\int_\M f (12 \cos^3\theta + 4 \cos^2\theta-6\cos\theta-2)d\theta.
\eqnn
In order to study the bracket $[v,\Delta ]$, we use two schemes $\SS$ and $\tilde{\SS}$ that discretize the diffusion-growth problem described above. 
We define $\tilde{\SS}$ similarly to $\SS$ (defined in Section \ref{s-ex}), but inverting steps 1 and 2. Hence $\SS$ does a series of growth and diffusion operations on the function $\mu_0$ starting with growth, while $\tilde{\SS}$ does the same starting with diffusion. Figure \ref{fig:Schemes} shows the first two iterations of each scheme, starting from the same function $\mu_0$ (renamed $x_0$ and $y_0$ for notation convenience), and denoting respectively by $x_n$ and $y_n$ the solutions after each iteration of $\SS$ and $\tilde{\SS}$.
We first apply this scheme to the initial signal given by the 
constant function $\mu_0(\theta)=0.1$.
%symmetric and smooth function $Mu_0(\theta)=\alpha (\cos\theta+1)$.
% According to \eqref{Brack0}, the initial bracket is given by: 
%\begin{equation}
%Mu_0 [\Delta, v] = - \beta \sin\theta \; \partial_\theta Mu_0 - 2\beta \cos\theta \; \partial^2_\thetaMu_0 = \beta\alpha ( \sin^2\theta + 2\cos^2\theta)=\beta\alpha ( 1 + \cos^2\theta). 
%\end{equation}
To numerically compute the lie bracket,
we apply the schemes $\SS$ and $\tilde{\SS}$ to $\mu_0$ and compute the numerical expression of the bracket given by $[\Delta,v]_{\text{num}} = \lim_{\epsilon\rightarrow 0} (y_1-x_1)/\epsilon^2$ (where $x_1$ and $y_1$ respectively correspond to the first iterations of $\SS$ and $\tilde{\SS}$).\\
%We visually compare the numerical and theoretical brackets, for values $\beta=1$, $\alpha=0.1$ (see Figure \ref{fig:Bracketnumtheo}).
%Figure \ref{fig:EvolutionConstantSignal} shows the evolution of the constant signal $\mu_0=0.1$ after one iteration of each scheme $\SS$ and $\tilde{\SS}$, with a vector field $v = (x-1,2y)$ and time-steps $t=\tau=0.05$ (so $T=t+\tau = 0.1$).
Figure \ref{fig:BracketnumConstantSignal} shows the convergence of the bracket when the time step $T=t+\tau$ tends to 0: The bracket converges to the theoretical value $[\Delta,v]_{\text{theo}}(\mu_0)=0.2\cos(\theta)$.\\
Figure \ref{fig:BracketnumCosineSignal} shows the convergence of the bracket for the initial signal $\mu_0(\theta)=0.1(\cos(\theta)+1)d\theta$. The bracket converges to the theoretical value $[\Delta,v]_{\text{theo}}(\mu_0)=12 \cos^3\theta + 4 \cos^2\theta-6\cos\theta-2$.
%\begin{figure}[h!]
%        \begin{center}
%                \includegraphics[trim=0cm 0cm 0cm 0cm, clip=true, scale=0.4]{EvolutionConstantSignal_b&w.eps}
%                        \caption{Evolution of the signal $\mu_0(\theta)=0.1$ after a time-step $T=0.1$ for the two schemes.}\label{fig:EvolutionConstantSignal}
%        \end{center}
%\end{figure}

\vspace{-1cm}

%\addtolength{\textheight}{-6cm}

\section*{CONCLUSION} 

In this paper, we described morphogenesis of organisms by evolution of a manifold, representing their boundary. The presence of morphogens, that are signals stimulating deformations, is modeled by the presence of a growth signal evolving on the manifold. \\
We showed that the resulting dynamics can be described by the reaction-diffusion Partial Differential Equation \r{e-main}, whose solution in the sense of measures is supported on the evolving manifold. 
The existence of solutions is achieved using the framework of Wasserstein distances.
In such equation, the interplay between the growth and diffusion can be described by a new concept of Lie bracket between the transport term and the Laplace-Beltrami diffusion operator. The two operators
are of different nature, but they can be both applied to measures on the ambient space.
This allows a precise definition of such Lie bracket and its explicit expression is given \r{e-brackint}. Then via numerical simulations we verified the effect of the non-commutativity of the diffusion and growth and the found expression for the bracket.\\
Future work will be devoted to develop a complete theory for PDEs on time-evolving manifolds
and numerically study the shapes achievable by appropriate control mechanisms.

\section*{ACKNOWLEDGMENT}

The authors acknowledge the partial support of the NSF Project "KI-Net", DMS Grant \# 1107444,
and of the National Institute of General Medical Sciences of the National Institutes of Health under Award Number R15GM101597, awarded to NY and BP.

%%%%%%%%%%%%%%%%%%%%%%%%%%%%%%%%%%%%%%%%%%%%%%%%%%%%%%%%%%%%%%%%%%%%%%%%%%%%%%%%

\end{document}